

\baselineskip=14pt
\parskip=10pt

\font\eightrm=cmr8 
\font\eighttt=cmtt8
\magnification=\magstephalf

\def\1{{\overline{1}}}
\def\2{{\overline{2}}}
\parindent=0pt
\overfullrule=0in

\def\frac#1#2{{#1 \over #2}}
\bf
\centerline
{
Teaching the Computer how to Discover(!) and then Prove(!!) (all by Itself(!!!)) 
}
\centerline
{
Analogs of Collatz's  Notorious 3x+1  Conjecture
}
\rm
\bigskip
\centerline{ {\it
Doron 
ZEILBERGER}\footnote{$^1$}
{\eightrm  \raggedright
Department of Mathematics, Rutgers University (New Brunswick),
Hill Center-Busch Campus, 110 Frelinghuysen Rd., Piscataway,
NJ 08854-8019, USA.
{\eighttt zeilberg  at math dot rutgers dot edu} ,
\hfill \break
{\eighttt http://www.math.rutgers.edu/\~{}zeilberg} .
First version: March 23, 2009.
Accompanied by Maple package {\eighttt LADAS}
downloadable  from
{\eighttt http://www.math.rutgers.edu/\~{}zeilberg/mamarim/mamarimhtml/collatz.html}, where one
can also find (very interesting!) output, consisting of 144 computer-generated theorems and proofs.
Supported in part by the NSF.
}
}

{\bf Mathematics: an Experimental Science}

In spite of their exponential growth, computer-assisted and especially computer-generated mathematical research
are still in their infancy. One approach is that of {\it formal proofs} using the {\it axiomatic method}.
This method is based on the {\it myth}, that goes back to Euclid,
that mathematics is a  {\it deductive} science, where one starts with 
a bunch of axioms, 
(initially supposed to be ``self-evident'', but later conceded as true-by-fiat) 
and then uses {\it rules of deduction},
and step-by-step,  arrives at (seemingly) non-trivial results.

Of course, mathematics {\it could} be presented that way, and unfortunately, often {\it is}.
But that is not how it is {\it discovered}. Pretending that the axiomatic method is how
mathematics should be done, and trying to indoctrinate poor computers to do it that way,
is a highly {\it inefficient} use of computers' time.

Deep inside, mathematics is, or at least should be, an {\it inductive} science. How do human mathematicians
come up with such amazing conjectures? By experimenting! How do they come up with such amazing
proofs? By experimenting!, and {\bf not} by combining axioms.

For many results in mathematics, one only needs one non-trivial ``axiom'', that of Peano's
``axiom'' of induction:
$$
P(0) \,\,\, \& \, \, \, ( P(n) \Rightarrow P(n+1)) \, \, \Rightarrow \, \,\, \, P(n) \quad for \quad all \quad n \geq 0 \quad .
$$
All we have to do then is to have the computer prove (all by itself) $P(0)$ (an a priori routine fact), and
the not-a-priori routine $P(n) \Rightarrow P(n+1)$, but, with the help of {\it symbolic computation},
treating $n$ as a {\it symbolic integer} (which is {\bf one} object) rather than as a {\it variable}
(ranging over positive integers), we might hope to have the computer do it for us.

Often one tries to prove a statement by induction, but {\it fails}. In that case one has
to {\it try again}. If $P(n)$ is inadequate to prove $P(n+1)$, perhaps we need
another statement $Q(n)$, such that $Q(n)$ and $P(n)$ imply $P(n+1)$. Alas, now we also
have to prove that $P(n)$ and $Q(n)$ imply $Q(n+1)$. This may be easier-said-than-done, and
we may be forced to introduce yet-another statement $R(n)$. If we want to avoid a
{\it Ponzi scheme}, we need this process to finally halt, but if we can train
our computer how to do it, it will be able to keep trying, before halting or giving up,
for a much  longer time.

Human mathematicians also do {\it symbol crunching}, most of which is purely routine, and
that is much better delegated to computers.

{\bf Research is one percent Inspiration and ninety nine percent Perspiration}

At this time of writing, humans are still needed to find {\it ideas} and {\it strategies}
for generating conjectures, and for proving these conjectures. But once the human
has some ideas, it is much more efficient to teach the computer these ideas, and let the
computer search for {\it conjectures}, and most impressively, {\it proofs}, all by itself!

In this case-study in computer-generated mathematics, I will 
describe how I read
{\it very carefully} a beautiful mathematical paper [AGKL], 
written by four {\it brilliant} human
beings: Ed Grove and Gerry Ladas, and their (at the time) 
respective students, Candy Kent
and Amal Amleh. I then extracted the {\it ideas}, 
looked at the {\it structure} of the  proofs,
but ignored the details. I then {\it taught} (programmed) my beloved 
computer to execute
these ideas, and since it is much more patient than a human being, it was able
to prove many more results. 
A very preliminary output can be viewed at the webpage of this
article

{\eighttt http://www.math.rutgers.edu/\~{}zeilberg/mamarim/mamarimhtml/collatz.html}  ,

where there are 144 computer-generated results, that include those proved in [AGKL].
The computer did {\it all} the phases of mathematical research all by {\bf itself}: 
conjecturing theorems, conjecturing proofs, and finally, verifying that the conjectured proofs are correct.

{\bf Background: The 3x+1 problem}

Mathematics abounds with easy-to-state yet (apparently) impossible-to-prove
statements, for example the Goldbach and twin-prime conjectures, but none
of them rivals the simplicity of statement, and most probably, difficulty
of proof, of the Collatz infamous 3x+1 conjecture, so beautifully exposited
in Jeff Lagarias' [L] masterpiece. I am sure that all the readers know it,
but for the sake of completeness, let me state it anyway.

{\bf The Collatz 3x+1 Problem}

Let $x(n), (n \geq 0$), be a sequence of positive integers defined by the 
{\it first-order} recurrence
$$
x(n+1) = \cases{
{{x(n)} \over {2}} ,& if \quad $x(n)$ \quad is \quad even ;\cr
              {{3x(n)+1} \over {2}} ,&  if \quad $x(n)$ \quad is \quad odd .\cr}   \quad,
$$
subject to the {\it initial condition} $x(0)=x_0$. Then
for {\it every} positive integer $x_0$, the sequence is eventually
periodic with the trivial cycle $(2,1)$.

Paul Erd\H{o}s claimed that mathematics is not yet ready for solving such
problems. I strongly believe that it very soon will be, thanks to the
emergence of computer-assisted and computer-generated methods, that
eventually would be able to {\it rigorously} prove such statements, 
for {\it all} $x_0 < \infty$,
as opposed to just numerically verifying it for $x_0<M$ for some large
specific integer $M$. Of course, we can't do it naively, since there
are infinitely many cases to check, and even our largest and fastest
computers are {\it finite} (and so are we, and our universe, for that matter).

If you try something and you fail, try and try again, then {\it generalize}
and/or {\it analogize}!, no use being a damn fool and trying to prove the
original conjecture. In 1995, Clark and Lewis[CL] played with
an analogous, {\it second-order} recurrence
$$
x(n+1) = \cases{{{x(n-1)+x(n)} \over {2}} ,
& if \quad $x(n-1)+x(n)$ \quad is \quad even ;\cr
             { -x(n-1)+x(n)},&  if \quad $x(n-1)+x(n)$ \quad is \quad odd .\cr}   \quad,
$$
and started experimenting with it, numerically,
by trying out some random initial values, and observing
what is going on.
For example, the following trajectory 
arises with $x(-1)=11,x(0)=16$:
$$
11, 16, 5, -11, -3, -7, -5, -6, -1, 5, 2, -3, -5, -4, 1, 5, 3, 4, 1, -3, -1, -2, -1, 1, 0, -1, -1, -1, -1, \dots  \quad  \quad ;
$$
while the following one, with $x(-1)=13$, $x(0)=6$:
$$
13, 6, -7, -13, -10, 3, 13, 8, -5, -13, -9, -11, -10, 1, 11, 6, -5, -11, -8, 3, 11, 7, 9, 8, -1, -9, -5, -7, -6,
$$
$$
1, 7, 4, -3, -7, -5, -6, -1, 5, 2, -3, -5, -4, 1, 5, 3, 4, 1, -3, -1, -2, -1, 1, 0, -1, -1, -1,-1,-1, \dots \quad  \quad ;
$$
or this one, with $x(-1)=11,x(0)=5$
$$
11, 5, 8, 3, -5, -1, -3, -2, 1, 3, 2, -1, -3, -2, 1, 3, 2, -1, -3 \dots \quad .
$$
They noticed that none of the trajectories they encountered seem to go to
infinity, but instead eventually end-up either being
identically $1$, or identically $-1$, or end-up with the six-period
$(-2,1,3,2,-1,-3)$ for ever after. 
They then (probably) wrote a short computer program
that checked their guess for all initial values 
$-100 \leq x(-1),x(0) \leq 100$, and they had ample {\it empirical}
evidence for the following

{\bf Theorem}([CL]) Let $x(-1)$ and $x(0)$ be any two integers
whose largest odd common divisor is $1$, then $x(n)$ is either
eventually the constant $1$, the constant $-1$, or the six cycle
$(-2,1,3,2,-1,-3)$.

They then worked fairly hard to give a purely human proof.

A few years later, in the remarkable paper [AGKL] already mentioned above,
Amleh, Grove, Kent and Ladas considered all $16$ recurrences,
for $\alpha, \beta, \gamma, \delta \in \{-1,1\}$. 
$$
x(n+1) = \cases{ {{\beta x(n-1)+ \alpha x(n)} \over {2}} ,
& if \quad $x(n-1)+x(n)$ \quad is \quad even ;\cr
             { \delta x(n-1)+ \gamma x(n)},&  
if \quad $x(n-1)+x(n)$ \quad is \quad odd .\cr}   \quad.
\eqno(AGKL)
$$

It is immediate to see that if the initial values $x(-1),x(0)$ have
an odd common factor, all the remaining terms would also, so
without loss of generality, one can consider the case
where $x(-1),x(0)$ do not have a common odd divisor. Also if they
have a common even divisor, as long as two consecutive terms
are even, the transformation is ``shrinking'', so without loss
of generality $gcd(x(-1),x(0))=1$.

By a lot of {\it human ingenuity} and a considerable
amount of {\it human toil}, Amleh et. al. proved that 
for two cases there are trajectories that go to infinity
($(\alpha, \beta, \gamma, \delta)=(1,1,1,1), (1,-1,1,-1)$),
and for ten other cases (out of the possible sixteen cases)
they proved that all trajectories end-up with
periodic orbits, and listed all of them. They left four open
cases. One case is
$$
(\alpha, \beta, \gamma, \delta)=(-1,1,1,1) \quad ,
$$
for which they conjectured that all trajectories end up
in one of the cycles
$$
(0,1,1) \quad, \quad (0,-1,-1) \quad , \quad
(2,5,7,1,-3,-2,-5,-7,-1,3) \quad ,
$$
and the case
$$
(\alpha, \beta, \gamma, \delta)=(-1,1,-1,-1) \quad ,
$$
for which they conjectured that all trajectories end up
in one of the cycles
$$
(-1,-2,3) \quad, \quad (1,2,-3) \quad , \quad
(-1,0,1,-1) \quad , \quad (1,0,-1,1) \quad .
$$
Later they found yet another cycle, of length $33$:
$$
(79, -31, -55, -12, 67, -55, -61, -3, 29, 16,  -45, 29, 37, 4, -41, 37, 39, 1, -19, -10, 29, -19, 
$$
$$
-24, 43, -19, -31, -6,  37, -31, -34, 65, -31, -48) \quad ,
$$
and made the modified conjecture that this completes the list.

The two other open cases are ``dual'' (see [AGKL]) to those two, and they would follow from them.

{\bf Deconstructing the Ingenious Human Proof}

For the sake of definiteness and exposition, let's focus on 
the original Clark-Lewis[CL] recurrence, 
$$
x(n+1) = \cases{{{x(n-1)+x(n)} \over {2}} ,
& if \quad $x(n-1)+x(n)$ \quad is \quad even ;\cr
             { -x(n-1)+x(n)},&  if \quad $x(n-1)+x(n)$ \quad is \quad odd .\cr}   \quad,
\eqno(CL)
$$
but following the proof strategy of [AGKL].

{\bf General Outline of the Amleh-Grove-Kent-Ladas Proof Strategy}

{\bf 1.} Try to prove that every trajectory is {\it bounded}.
It would then follow immediately from the pigeon-hole principle,
that there exists a pair of integers that shows up twice
as $(x(n),x(n+1))$, and once that happens it is trapped for ever
into an orbit, but {\it a priori} infinitely many ultimate periods may show up,
for different pairs of initial values.

{\bf 2.} In order to prove the boundedness assertion of {\bf 1} more
specifically, look for two integers, $c_1$ and $c_2$ such that,
empirically for now:
$$
|x(n)| \leq c_1 |x(-1)|+ c_2 |x(0)| \quad ,
\eqno(FundIneq)
$$
(for the present case $c_1=1,c_2=1$ works).

{\bf 3.} Try to prove $(FundIneq)$ by induction. You will
probably fail at first, but whenever you fail, you can add more
hypotheses, and try to prove a more general statement, proving
several inequalities at once by induction, including the one
you really need. You may have to keep adding
more and more statements, and this process might never end
(and indeed if $(FunIneq)$ is false, it better not end!), but
if in luck, this would end. We will describe it in more detail later.

{\bf 4.} We now know that every trajectory must end-up with a cycle. We want to
characterize {\it all} these cycles. Since a cycle is ``cyclic'', we can, without
loss of generality, make the {\bf largest} element (in absolute value) be the
$x(1)$ . Also since all the states are mod $2$, we can
take the initial conditions
$x(-1)=a$ and $x(0)=b$ or $x(-1)=a$ and $x(0)=-b$ , where $a$ and $b$ denote (symbolic) positive integers.
(We should separately treat the cases $b=0$ and $a=0$, but this is really easy).

By cyclicity (considering the last two elements $x(n-1),x(n)$ as our $x(-1),x(0)$, and $x(1)$ as an entry that
comes after), we have the {\bf necessary condition}
$$
|x(1)| \leq c_1 |x(n-1)|+ c_2 |x(n)| \quad ,
$$
that in addition to $|x(n-1)| \leq |x(1)|$, $|x(n)| \leq |x(1)|$, 
enable us to ``rule out'' lots of possible parity sequences defined below.

{\bf 5.} 
Define the {\bf parity sequence} of a trajectory $x(n)$ to
be the sequence $x(n)\, mod \, 2$. For example, the 
parity sequence of the trajectory
$$
11, 5, 8, 3, -5, -1, -3, -2, 1, 3, 2, -1, -3, -2, 1, 3, 2, -1, -3, \dots 
$$
is
$$
1, 1, 0, 1, 1, 1, 1, 0, 1, 1, 0, 1, 1, 0, 1, 1, 0, 1, 1, \dots  \quad .
$$
We are looking for candidate parity sequences that {\bf do not} violate the
{\bf necessary condition} of {\bf 4}. Note that if we know the parity sequence
of a trajectory, we can express each term of the trajectory, in particular the
last two, as certain explicit linear combinations of $a=x(-1), b=x(0)$.

{\bf 6.} Using {\bf 4}, we first empirically discover ``all'' the
possible parity sequences of potential cycles.
These parity sequences turn out (in all the 
successful cases encountered so far) to be {\it regular expressions} parametrized
by few integer parameters. For example, for the present case they turn out to be
$$
(011)^{m_1} 1^{m_2} \quad, \quad (1011) 1^{m_1} \quad , \quad
(011)^{m_1} 1^{m_2} 0 \quad , \quad (011)^{m_1} 1^{m_2} 01 \quad.
$$

{\bf 7.} Rigorously prove, by induction on $m_1,m_2, \dots$, that these are the only ones,
by using the inequalities of {\bf 4}, but with symbolic $m_1, m_2, \dots$
(in addition to the symbolic $a,b$).

{\bf 8.} You can now express in {\bf closed form}, in terms
of $a:=x(-1)$ and $b:=x(0)$, the general two last terms,
let's call them $x(n-1),x(n)$, of such
a trajectory. To investigate whether it can form a periodic orbit,
try to solve the system of two linear equations
$$
x(n-1)=a \quad , \quad x(n)=b \quad , \quad (or x(n-1)=a \quad , \quad x(n)=-b) \quad ,
$$
getting a homogeneous set of two equations with two unknowns.
$a$ and $b$. Set the determinant of the $2 \times 2$ matrix of the
resulting system of equations 
to zero, and get an expression in $m_1,m_2, \dots$ that is usually (i.e. generically) not zero.
For those rare cases where it is $0$, you would get
a hope for a cycle, and solving for $a$ and $b$ would give them to you.
(This is a bit analogous to finding eigenvalues, followed by the corresponding eigenvectors).

Before we can teach the computer how to discover and prove Collatz-type theorems, we need
to teach ourselves, in somewhat greater detail,  the brilliant ideas of [AGKL].Let's do this step-by-step for
the original Clark-Lewis recurrence
$$
x(n+1) = \cases{{{x(n-1)+x(n)} \over {2}} ,
& if \quad $x(n-1)+x(n)$ \quad is \quad even ;\cr
             { -x(n-1)+x(n)},&  if \quad $x(n-1)+x(n)$ \quad is \quad odd .\cr}   \quad,
\eqno(CL)
$$

{\bf Proving Boundedness}

By using procedure {\tt FindPreScheme} in the Maple package {\tt LADAS}, the computer quickly makes the
conjecture
$$
|x(n)| \leq |x(-1)|+ |x(0)| \quad , \quad for \quad n \geq -1 \quad .
$$
The natural approach would be to use {\it induction} on $n$. But it turns out that we (usually) need a {\bf stronger}
statement, so that we will have more {\it elbow-room} for the inductive argument.

We will shortly describe how the computer can {\it dynamically} construct this stronger statement,
but let us first describe how the computer (and of course, us) can rigorously prove its validity,
once it has been formulated.
We will give more details than in customary in human discourse, avoiding the phrases and words ``it is easily seen'',
``obviously'' (and whenever possible) ``similarly''. This level of detail is necessary if we want to be able
to transfer our know-how to our silicon brethern.

In the present case (recurrence $(CL)$), the stronger statement turns out to be as follows.
Let $A=|x(-1)|+|x(0)|$.

Case I(n): If $x(n-1),x(n)$ are both {\it odd} then: 
$$
|x(n-1)| \leq A \quad , \quad |x(n)|  \leq A \quad .
$$
Case II(n): If $x(n-1)$ is {\it odd} and $x(n)$ is  {\it even}, then:
$$
|x(n-1)| \leq A \quad , \quad  |x(n)|  \leq A \quad , \quad |-x(n-1)+x(n)| \leq A \quad .
$$
Case III(n): If $x(n-1)$ is {\it even} and ,$x(n)$ is  {\it odd},  then
$$
|x(n-1)| \leq A \quad , \quad |x(n)|  \leq A \quad , \quad |-x(n-1)+x(n)| \leq A \quad .
$$

(In the present example of $(CL)$, the two last cases can be combined into ``if $x(n-1)+x(n)$ is odd'', since the three inequalities
are identical. Usually this is not the case).

The natural approach would be to prove it by induction. There are three {\bf states}

$I(n)$: $(x(n-1),x(n))$=(odd, odd) \quad ,

$II(n)$: $(x(n-1),x(n))$=(odd, even), 

$III(n)$: $(x(n-1),x(n))$=(even,odd).

Let's consider them each at their turn.

If we are currently in state $I(n)$
(i.e. $(x(n-1),x(n))$=(odd, odd)), then the next ``state'' is either $I(n+1)$ (i.e. $(x(n),x(n+1))$=(odd, odd)) or
$II(n+1)$ (i.e. $(x(n),x(n+1))$=(odd, even)).

Case Ia: $(x(n-1),x(n))$=(odd, odd) \quad , \quad $(x(n),x(n+1))$=(odd, odd) .

We have to prove statement $I(n+1)$, in other words:  $|x(n)| \leq A \quad , \quad |x(n+1)|  \leq A$ .

In terms of $x(n-1),x(n)$ we have to prove: $|x(n)| \leq A \quad , \quad |{{1} \over {2}} x(n-1)+ {{1} \over {2}} x(n) |  \leq A$ .

The first inequality is already contained in the (inductive) premise, while for the second we use the
{\bf triangle inequality}
$$
|{{1} \over {2}} x(n-1)+ {{1} \over {2}} x(n) | \leq {{1} \over {2}} |x(n-1)|+ {{1} \over {2}} |x(n)| \quad,
$$
By the {\it inductive hypothesis} we know that $|x(n-1)| \leq A \quad$ , $\quad |x(n)| \leq A$, so it follws that
$$
|{{1} \over {2}} x(n-1)+ {{1} \over {2}} x(n)| 
\, \, \leq  \,\, {{1} \over {2}} |x(n-1)|+ {{1} \over {2}} |x(n)| \,\, \leq \,\, {{1} \over {2}}A+ {{1} \over {2}} A \leq A \quad ,
$$
since ${{1} \over {2}}+{{1} \over {2}} \leq 1$.

Case Ib: $(x(n-1),x(n))$=(odd, odd) \quad , \quad $(x(n),x(n+1))$=(odd, even).

We have to prove $II(n+1)$, in other words:
$$
|x(n)| \leq A \quad , \quad |x(n+1)|  \leq A \quad , \quad |-x(n)+x(n+1)| \leq A \quad .
$$
Expressed in terms of $x(n-1)$ and $x(n)$ these are
$$
|x(n)| \leq A \quad , \quad |(x(n-1)+x(n))/2|  \leq A \quad , \quad |-x(n)+(x(n-1)+x(n))/2| \leq A \quad .
$$
Cleaning up, we have to prove:
$$
|x(n)| \leq A \quad , \quad |{{1} \over {2}}x(n-1)+{{1} \over {2}}x(n)|  \leq A \quad , \quad |{{1} \over {2}} x(n-1)-{{1} \over {2}}x(n)| \leq A \quad .
$$
The first inequality is contained in the inductive hypothesis, the second one is identical to the one done above,
while, for the third one, we have, once again by the triangle inequality:
$$
|{{1} \over {2}} x(n-1)- {{1} \over {2}} x(n)| \,\, \leq  \,\,|{{1} \over {2}}| A+ |-{{1} \over {2}}| A \,\, \leq \,\, A \quad ,
$$
since $|{{1} \over {2}}|+ |{{-1} \over {2}}| \leq 1$.

If we are currently in state II(n), the next state is necessarily III(n+1). We  have to prove 

III(n+1): $|x(n)| \leq A, |x(n+1)|  \leq A, |-x(n)+x(n+1)| \leq A$.

In terms of $x(n-1),x(n)$ these are:  $|x(n)| \leq A, |-x(n-1)+x(n)|  \leq A, |-x(n)-x(n-1)+x(n)| \leq A$.

Simplifying, these are: $|x(n)| \leq A, |-x(n-1)+x(n)|  \leq A, |-x(n-1)| \leq A$. 

The first two inequalities are
part of the inductive hypothesis, while the third one follows from the deep fact that $|-1|=1$:
$$
|-x(n-1)|=|(-1)x(n-1)|= |-1||x(n-1)|= 1 |x(n-1)|=|x(n-1)| \leq A \quad .
$$

If we are currently in state III(n), the next state is necessarily II(n+1). In general, this would be different than
the previous case, but in this example it is identical, so I hope that the reader will forgive me for using ``similarly'',
since this is not only {\it similarly} (in the colloquial human sense) but {\it identical}, and  even a computer
can realize that, and prevent a duplication  of effort.

{\bf How to teach the Computer to Argue as above}

In the above proof, we needed to know 

(i) If we are currently at a certain $state(n)$, what state(s) is (are) next?

(ii) Rewrite the expressions given in terms of $x(n),x(n+1)$, in terms of $x(n-1),x(n)$, using the rules of the
recursion, according to the current state.

(iii) Know how to apply the triangle inequality, or realize that
the inequality that we have to prove is already part of the inductive hypothesis. 

Tasks (i) and (ii) are obvious, while for (iii) we need to prove that a bunch of inequalities
$$
|L_{1}| \leq A \quad , \quad  |L_{2}| \leq A \quad , \quad \dots \quad , \quad  |L_{m}| \leq A
$$
imply another one
$$
|M| \leq A \quad .
$$
If in luck, $M$ is already one of the $L_i$'s, or $M=-L_i$, for some $i$. Otherwise, we do a double do-loop looking
for $i$ and $j$ such that one can write
$$
M=c_1 L_i + c_2 L_j \quad ,
$$
and then we check that $|c_1|+|c_2| \leq 1$. If this fails, we try a triple do-loop, looking for $1 \leq i<j<k \leq n$ such that
$$
M=c_1 L_i + c_2 L_j +c_3 L_k \quad ,
$$
$|c_1|+|c_2| +|c_3| \leq 1$, etc.

All this involves solving simple systems of linear equation, that Maple can do very fast.

{\bf OK, A Computer can automatically prove, by Induction, such a ``Scheme'' of inequalities, but we sure
need humans to come up with it!}

{\bf Wrong, of course!} The computer starts with just what it really wants to know, namely
$|x(n)| \leq A$, and tries to prove it by induction. Chances are that there wouldn't be
enough assumptions, so to paraphrase Guru Greg Chaitin, we ``add it as an axiom'', i.e.
hypothesis, and keep trying to prove our ``partial scheme'' by induction, and whenever we need
another assumption, let's just add it. If we are lucky, this process would eventually halt,
and we would be done.

Let's illustrate it with the $(CL)$ recurrence. We start with the single inequality $|x(n)| \leq A$ in all states.
Obviously we need to add $|x(n-1)| \leq A$ right away. This suffices for the transition
$I(n) \Rightarrow I(n+1)$, and  $I(n) \Rightarrow II(n+1)$. Alas for
$II(n) \Rightarrow III(n+1)$ and $III(n) \Rightarrow II(n+1)$ we can't prove $|-x(n-1)+x(n)| \leq A$,
so we {\it add} it as part of the statement, and as we saw above, this suffices.

In general, we may have to add more and more hypotheses, but if all goes well, we would eventually converge.

{\bf Conjecturing the possible parity sequences of Cycles}

Procedure {\tt DAPf} of our Maple package {\tt LADAS} constructs all feasible parity sequences
of a given rule, of any given length. This is a recursive procedure that starts out with
all feasible symbolic trajectories of length 3, phrased in terms of $a$ and $b$ that conform to
the convention that the third element, $x(1)$, is the largest in absolute value, and using $a$ and $b$ to denote
{\it symbolic} {\bf positive} integers. For example for Rule 6 of [AGKL]
$$
x(n+1) = \cases{{{-x(n-1)+x(n)} \over {2}} ,
& if \quad $x(n-1)+x(n)$ \quad is \quad even ;\cr
             { -x(n-1)+x(n)},&  if \quad $x(n-1)+x(n)$ \quad is \quad odd .\cr}   \quad,
\eqno(E6)
$$
the feasible trajectories (with the above convention that the third entry, $x(1)$ is largest in
absolute value) are:
$$
[[a, -b, -a - b], [1, 2, 1]], [[a, -b, -a - b], [2, 1, 1]] \quad,
$$
(here we use $2$ instead of $0$, since in the program $0$ denotes ``either odd or even'').
So we already know that it can't start with both entries being odd (since then $x(1)$ wouldn't be
able to be larger than both $a$ and $b$). How does our computer generate the possible trajectories
of length $4$?.
For each parity sequence of length $3$ it considers the two 
extensions obtained by appending $2$ (even) and $1$(odd). Then Maple {\it automatically} solves the system
of linear inequalities (featuring absolute values, Maple can do it!) and if there are no solutions,
discards that option. Step by step it constructs the potential parity sequences,
and also expresses all the entries symbolically.

After  Maple gathered enough data, it looks for {\it patterns}, and 
by using procedure {\tt gGREplus} (of the Maple package {\tt LADAS}),
the computer finds, empirically, all the {\it regular expressions} in
the alphabet $\{0,1\}$ that are parity sequences of trajectories that do not violate
any of the necessary conditions of {\bf 4}.
To {\it prove} that these are the
{\it only} possibilities, rigorously, the computer first (symbolically) computes,  
(by using powers of matrices, or guessing, using the obvious format and then proving it by induction),
{\it explicit expressions} for the last two elements of a trajectory $x(n-1),x(n)$.
For example, for $(E6)$ (resorting back to the convention of the paper, rather than the program, of denoting 
an even integer by $0$ rather than $2$)
there are three such possible symbolic regular expressions:
$$
(011)^{m_1} \quad, \quad (011)^{m_1}1 \quad, \quad (011)^{m_1}01 \quad .
$$
To prove that these are indeed all of them, simply use induction combined with the inequality of {\bf 4}
and rule out, for example the possibility $(011)^{m_1}0$, that would violate (for symbolic $m_1$ as well
as symbolic $a$ and $b$) that inequality. This too {\it can} be done automatically.

{\bf Finding All the Cycles}

Now that we know what symbolic parity sequences may show up, we can
explicitly express the last two terms of a trajectory (that is a potential cycle) 
in terms of $a$, $b$ and $m_1, m_2 , \dots $. Simply use symbolic powers of
matrices (that Maple can easily do by finding the eigenvalues of the relevant
matrices). For example, for rule $(E6)$ and the symbolic
trajectory $(011)^{m_1}$, we have
$$
x(-1)=a \quad , \quad x(0)=-b \quad , \quad x(n-1)=(-1)^{m_1}+(-1)^{m_1+1} b \quad , \quad x(n)=a(-{{1} \over {2}})^{m_1} \quad .
$$
To investigate whether there exist $m_1$ (a positive integer) and positive integers $a$ and $b$ such that
this forms a cycle, we have to solve:
$$
x(-1)=x(n-1) \quad , \quad x(0)=x(n) \quad .
$$
In other words
$$
a=a(-1)^{m_1}+(-1)^{m_1+1} b \quad , -b=a(-{{1} \over {2}})^{m_1} \quad .
$$
Moving everything to the left:
$$
a(1-(-1)^{m_1})+(-1)^{m_1} b=0 \quad, \quad a(-{{1} \over {2}})^{m_1}+b=0 \quad .
$$
In order for there to be a non-zero solution $(a,b)$, the determinant:
$$
1-(-1)^{m_1}+{{1} \over{2^{m_1}}} =0 \quad,
$$
must vanish, but that can never happen. Similarly the other two (potential)
symbolic trajectories can be ruled out. The only contributors are those where
either $a=0$ or $b=0$, and these easy cases are treated separately.

{\bf Output}

My beloved computer, Shalosh B. Ekhad, used my Maple package {\tt LADAS} to
consider, more generally, {\it all} $256$ possible second-order difference equations of the form
$$
x(n+1) = \cases{ 
{{\alpha_1 x(n-1)+ \alpha_2 x(n)} \over {2}} , & if \quad $x(n-1)$ \quad is \quad even  \quad and \quad $x(n)$ is even  ;\cr
{{\alpha_3 x(n-1)+ \alpha_4 x(n)} \over {2}} , & if \quad $x(n-1)$ \quad is \quad odd  \quad and \quad $x(n)$ is odd  ;\cr
{ \alpha_5 x(n-1)+ \alpha_6 x(n)},&   if \quad $x(n-1)$ \quad is \quad odd \quad and \quad $x(n)$ \quad is \quad even.\cr
{ \alpha_7 x(n-1)+ \alpha_8 x(n)},&   if \quad $x(n-1)$ \quad is \quad even \quad and \quad $x(n)$ \quad is \quad odd.\cr
}   \quad,
$$
for $(\alpha_1, \dots , \alpha_8) \in \{-1,1\}^8$.
It was successful in $144$ cases (more than a half!). All these computer-generated theorems (and proofs!),
complete with the proving schemes for the fundamental inequalities (and their detailed proofs!), can be viewed from the webpage of this
article

{\tt http://www.math.rutgers.edu/\~{}zeilberg/mamarim/mamarimhtml/collatz.html} .

{\bf The Maple package COLLATZ} 

 In addition to the main package {\tt LADAS}, we also developed
a purely empirical Maple package, {\tt COLLATZ}, to conjecture the cycles of generalized Collatz-type transformations
of the form $n \rightarrow a_i n + b_i$ if $n \equiv \, i \,\, mod \,\, m$, ($i=1 \dots m)$ such that
$ma_i$ and $ia_i+b_i$ are integers.

{\bf Further Work} 

This is but the {\it tip of an iceberg}, except that the specific weight of 
{\it our} ice is  much closer to $1$ than that of real ice. Obvious extensions would be to
consider general parity-dependent transformations on $Z^2$, rather than
those of the form $(a,b) \rightarrow (b, F(a,b))$ considered here.
More generally, we should consider transformations of the form
$$
(a,b) \rightarrow F(a,b) \quad,
$$
where $F(a,b)=F_{i,j}(a,b)$ if $a \equiv i \,\,\, mod \,\,\, m$, $b \equiv j \,\,\, mod \,\,\, m$,
and we have $m^2$ linear transformations that map pairs $(a,b)$ such that
if $a$ mod $m$ and $b$ mod $m$ equal $i$ and $j$ respectively, then $F_{i,j}(a,b)$ outputs
a pair of integers.

Why stop at transformations of $Z^2$? We can also consider transformations of $Z^d$ for
$d=3,4, \dots$.

The two remaining open cases of [AGKL], that we can't do either, may very well be amenable to
the present approach, except that the $c_1,c_2$ that feature in the
inequality $|x(n)| \leq c_1 |x(-1)|+ c_2 |x(0)|$ are rather large, and our program was not
fast enough to conjecture them? It would be a good idea to take our amateurish and slow
program and speed it up (say using Java and C), and perhaps be able to tackle these.
On the other hand, we may need yet another idea, and perhaps the two remaining cases of [AGKL] are
just as hard as the original Collatz problem.

What about {\it affine-linear} transformations? The original Collatz recurrence is a very simple
first-order recurrence given in that way. What about using the [AGKL] method,
as computerized here, to systematically explore even one-dimensional analogs of the $3x+1$ problem,
given by higher moduli, rather than two? See the Maple package {\tt COLLATZ} mentioned above,
for very preliminary (empirical) investigations.

And who knows? Perhaps a sufficiently higher-order recurrence, amenable to the present approach
(or to a yet-to-be-discovered extension) would imply, via an appropriate specialization,
the good-old Collatz problem? 

So there is a lot to do. Research is a relay race. Myself (and even my beloved Shalosh) are already tired
of this project,
but we do hope that other people (and machines!) will take over. Amen.

{\bf References}

[AGKL] A.M. Amleh, E.A. Grove. C.M. Kent, and G. Ladas, {\it On some difference equations with eventually 
periodic solutions}, J. Math. Anal. Appl. {\bf 223}(1998), 196-215.

[CL] D. Clark and J.T. Lewis, {\it A Collatz-type difference equation}, Congr. Numer. {\bf 111}(1995), 129-135.

[L] J.C. Lagarias, {\it The $3x+1$ problem and its generalization}, Amer. Math. Monthly {\bf 92} (1985), 3-23.

\end